\documentclass[runningheads]{llncs}
\usepackage{amsmath}
\usepackage{amssymb}  

\usepackage{graphicx}
\usepackage{xcolor}

\usepackage{mathtools}
\usepackage{caption}
\usepackage{subcaption}

\numberwithin{equation}{section}
\numberwithin{figure}{section}

\begin{document}
\title{Discrete gradient methods \\for irreversible port-Hamiltonian systems}

\author{Alexandre Anahory Simoes\inst{1} \and
David Martín de Diego\inst{2*} \and \\
Bernhard Maschke\inst{3}}
\authorrunning{A. Anahory Simoes et al.}
%
\institute{IE School of Sciences and Technology, P.º Castellana 259, 28029, Madrid, Spain \email{alexandre.anahory@ie.edu } \and
Instituto de Ciencias Matematicas (CSIC-UAM-UC3M-UCM), Calle Nicolas Cabrera 13-15, 28049 Madrid, Spain
\email{david.martin@icmat.es}\\ \and
 LAGEPP, UCB Lyon 1 - CNRS UMR 5007,
CPE Lyon - Bâtiment 308 G, Université Claude Bernard Lyon-1, 43,
Bd. du 11 Novembre 1918, F-69622 Villeurbanne, France
\email{bernhard.maschke@univ-lyon1.fr}
}
\maketitle              
\begin{abstract}

In this paper we introduce discrete gradient methods to discretize  irreversible port-Hamiltonian systems showing that the main qualitative properties of the continuous system  are preserved using this kind discretizations methods. 

\end{abstract}
\section{Introduction}

In \cite{Ramirez_ChemEngSci13,Ramirez_EJC13,RamirezThesis:2012} irreversible
port-Hamiltonian Systems have been suggested as a class of port-Hamiltonian systems \cite{Geoplex09} adapted to represent open physical systems in a way compatible with irreversible Thermodynamics \cite[chap.14]{Callen85}.
This formulation uses two thermodynamic potentials, the total energy and the entropy functions, but is formulated as a Hamiltonian system defined on an almost-Poisson bracket, contrary to other formulations which use a metriplectic structure \cite{Morrison86,Grmela_JPhysComm_18}.
Furthermore is encompasses open thermodynamic systems by defining
port variables and port maps associated with reversible and irreversible
phenomena at the interface of the system \cite{Kirchhoff_submIFAC_WC22b}.
In this paper, we shall define geometric integrators of these irreversible port-Hamiltonian systems which preserve the
energy and entropy balance equations
as it has been done for contact Hamiltonian systems  \cite{ALLM2020} or metriplectic systems \cite{Oettinger_JNT_18}.
 We shall use discrete gradient methods that are typically designed for numerical integration   of ordinary differential equations (ODEs), while ensuring the preservation of  certain structures of the continuous system as for instance energy preservation, dissipation, as it has also been used for Port Hamiltonian Systems \cite{Falaize_ApplSci_2016_PassSimAudio,Kotyczka_SCL_19_DiscreteTimePHS}. We will show, that the suggested method preserves both the invariance of energy and the irreversibility property of original system (see \cite{GONZ,ITOH,McQuRo} and references therein).

\section{Irreversible Port-Hamiltonian   systems}
\label{section:iphs} Irreversible Port Hamiltonian Systems have been
introduced in \cite{Ramirez_ChemEngSci13,Ramirez_EJC13}, as a class
of control systems for irreversible thermodynamic systems. Their
Hamiltonian function is the total energy of the system and they are
defined with respect to an almost-Poisson bracket which allows to
embed the second principle of Irreversible Thermodynamics, namely
the irreversible entropy creation term due to the irreversible phenomena.
Furthermore they represent also, open systems by expressing the interactions
of the system with its environment by using pairs of conjugated port
variables, called port variables. The definition given below, considers
a more precise expression, introduced in \cite{Kirchhoff_submIFAC_WC22b},
than in \cite{Ramirez_ChemEngSci13,Ramirez_EJC13}, of the input map
associated with a irreversible phenomena occurring at the interface
of the system such as the heat conduction.

\subsection{Definition}
\begin{definition}
\label{def:IPHS-IrrPortMap} An \textbf{irreversible port-Hamiltonian
Systems (IPHS)} is the nonlinear control system 
\begin{align}
\left(\begin{array}{c}
\frac{dx}{dt}\\
-y
\end{array}\right)= & \left[\sum_{i=1}^{k}\gamma_{i}\left(x,\tfrac{\partial H}{\partial x}\right)\left\{ S,H\right\} _{J_{i}}\left(\begin{array}{cc}
J_{i} & 0\\
0 & 0
\end{array}\right)+\sum_{\alpha=1}^{l}\left(\begin{array}{cc}
M_{\alpha} & 0\\
0 & 0
\end{array}\right)\right.\label{eq:IPHS_Matrix}\\
 & +\sum_{j=1}^{\bar{k}}\gamma_{\textrm{port},j}\left(x,\tfrac{\partial H}{\partial x},u\right)\left\{ S_{\textrm{tot}},H_{\textrm{tot}}\right\} _{J_{\textrm{port},g_{j}}}\left(\begin{array}{cc}
0 & g_{j}\\
-g_{j}^{\top} & 0
\end{array}\right)\nonumber \\
 & \left.+\sum_{\beta=1}^{\bar{l}}\left(\begin{array}{cc}
0 & g_{S\beta}\\
-g_{S\beta}^{\top} & 0
\end{array}\right)\right]\left(\begin{array}{c}
\frac{\partial H}{\partial x}\\
u
\end{array}\right)\nonumber 
\end{align}
where $x\left(t\right)\in\mathbb{R}^{n}$ is the \emph{state vector},
$u\left(t\right)\in\mathbb{R}^{m}$ is the \emph{control input} vector,
and defined by 
\begin{enumerate}
\item[(i)] two (smooth) real functions called \emph{Hamiltonian function} $H\in C^{\infty}(\mathbb{R}^{n})$
and \emph{entropy function} $S\in C^{\infty}(\mathbb{R}^{n})$
\item[(i)] two sets of structure matrices $J_{i}(x)\in\mathbb{R}^{n\times n}$,
$1\leq i\leq k$ and $M_{\alpha}(x)\in\mathbb{R}^{n\times n}$, $1\leq\alpha\leq m$,
defining \emph{almost-Poisson brackets} $\{\;,\;\}_{J_{i}}$ and $\{\;,\;\}_{M_{\alpha}}$,
respectively, and the \emph{real, strictly positive functions} $\gamma_{i}(x,\frac{\partial H}{\partial x})=\hat{\gamma}_{i}(x)\in C^{\infty}(\mathbb{R}^{n},\,\mathbb{R}_{+}^{*})$
.
\item[(iii)] the irreversible port-maps are defined by the \emph{input matrices}
$g_{j}(x)$ , $1\leq j\leq\bar{k}$, and the real strictly positive
functions. $\gamma_{\textrm{port},j}\left(x,\frac{\partial H}{\partial x},u\right)=\widehat{\gamma_{\textrm{port},j}}(x,u)\in C^{\infty}(\mathbb{R}^{n+m},\,\mathbb{R}_{+}^{*})$
. The \emph{port bracket} is defined with respect to the additional
functions $H_{\textrm{tot}}\left(x,\xi\right)=H\left(x\right)+u^{\top}\xi\;,\xi\in\mathbb{R}^{m}$
and $S_{\textrm{tot}}\left(x,\xi\right)=S\left(x\right)+\tau^{\top}\xi\;,\tau\in\mathbb{R}^{m}\;$
as 
\[
\left\{ S_{\textrm{tot}},H_{\textrm{tot}}\right\} _{J_{\textrm{port},g}}=\left[\left(g^{\top}\frac{\partial S}{\partial x}\right)^{\top}u-\tau^{\top}\left(g^{\top}\frac{\partial H}{\partial x}\right)\right]
\]
\item[(iii)] the reversible port-maps defined by the \emph{input matrices} $g_{S\beta}(x)$
in $\mathbb{R}^{n\times m}$, $1\leq\beta\leq\bar{l}$ . We assume
additionally that $\left\{ S_{\textrm{tot}},H_{\textrm{tot}}\right\} _{J_{\textrm{port},g_{S\beta}}}=0$
for all $1\leq\beta\leq\bar{l}$ 
\end{enumerate}
\end{definition}

\begin{remark}
Note that the definition given here is more general than the one given
in \cite{Kirchhoff_submIFAC_WC22b}, in the sense that the structure
matrices of the port maps may be dependent on the state.
\end{remark}

From the definition of the Irreversible Port Hamiltonian Systems,
one deduces the following energy and entropy balance equations 
\begin{align*}
 & \frac{dH}{dt}-y^{\top}u=0\\
 & \frac{dS}{dt}-y^{\top}\tau=\sum_{i}\hat{\gamma_{i}}\left(x\right)\left\{ S,H\right\} _{J_{i}}^{2}+\sum_{j}\widehat{\gamma_{\textrm{port},j}}\left(x,u\right)\left\{ S_{\textrm{tot}},H_{\textrm{tot}}\right\} _{J_{\textrm{port},g_{j}}}^{2}\geq0
\end{align*}
 where the right-hand side of the second equation is equal to the
irreversible entropy creation due to the irreversible phenomena in
the system and at its boundary.

\subsection{Example: gas-piston system}

This system describes the thermodynamic properties of a gas contained
in a cylinder closed by a piston submitted to gravity (\cite{Ramirez_EJC13})
and an an external force $u_{2}$, the gas is subject to an exchange
of heat with a thermostat at temperature $u_{1}$, obeying Fourier's
heat conduction law with coefficient $\lambda_{e}>0$. Mathematically
this system is described giving the mechanical Hamiltonian $H_{mec}(q,p)=\frac{1}{2m}p^{2}+mgq\;.$
where $q$ describes the altitude of the piston and $p$ its associated
momentum. On the other hand, the properties of the perfect gas may
be defined by its internal energy $U(S,V)$ where $S$ is the entropy
variable and $V$ the volume variable. Therefore, the total energy
if 
\begin{equation}
H(S,V,q,p)=U(S,V)+H_{mec}(q,p)\label{eq:HamiltonianGasPiston}
\end{equation}
Denote by $\frac{\partial H}{\partial S}=T,\quad\frac{\partial H}{\partial V}=-P,\quad\frac{\partial H}{\partial q}=mg,\quad\frac{\partial H}{\partial p}=v,\quad$
the corresponding co-energy variables, where $T$ represents the temperature
of the gas and $P$ its pressure, $mg$ is the gravity force and $v$
is the velocity of the piston. We assume that when the piston moves,
a non-adiabatic transformation due to mechanical friction and the
viscosity of the gas, induces a heat flow in the gas only, leading
to the entropy balance law $\frac{dS}{dt}=\frac{1}{T}\mu v^{2}\geq0$
parametrized by $\mu\geq0$.

The equations of motion of the gas-piston system can be written as:
\[
\left(\begin{array}{r}
dS/dt\\
dV/dt\\
dq/dt\\
dp/dt
\end{array}\right)=\left(\begin{array}{cccc}
0 & 0 & 0 & \mu v/T\\
0 & 0 & 0 & A\\
0 & 0 & 0 & 1\\
-\mu v/T & -A & -1 & 0
\end{array}\right)\left(\begin{array}{r}
T\\
-P\\
mg\\
v
\end{array}\right)+\left(\begin{array}{c}
\lambda_{e}\left(\frac{1}{T}-\frac{1}{u_{1}}\right)\\
0\\
0\\
u_{2}
\end{array}\right)
\]
where $A$ denotes the area of the piston. This system is an Irreversible
Port-Hamiltonian Systems generated by the Hamiltonian (\ref{eq:HamiltonianGasPiston})
and the entropy function, denoted with an abuse of notation $S(S,V,q,p)=S$.
The reversible coupling between the kinetic and potential mechanical
energy as well as with the gas through the piston are given by the
skew-symmetric matrix $M$ and the coupling between the mechanical
displacement and the heat flow due to the dissipation is given by
the skew-symmetric matrix $J$ defined by 
$$
J=\left(\begin{array}{cccc}
0 & 0 & 0 & 1\\
0 & 0 & 0 & 0\\
0 & 0 & 0 & 0\\
-1 & 0 & 0 & 0
\end{array}\right),\ M=\left(\begin{array}{cccc}
0 & 0 & 0 & 0\\
0 & 0 & 0 & A\\
0 & 0 & 0 & 1\\
0 & -A & -1 & 0
\end{array}\right)
$$
One may check that $\nabla S=(1,0,0,0)\in\ker M$, that is, the entropy
function is a Casimir function for $M$. Furthermore the Poisson backet
is $\{S,H\}_{J}=v$ which is the velocity of the piston and is the
driving force of the dissipation. The positive function giving the
constitutive relation of the dissipation is $\hat{\gamma}(S,V,q,p)=\frac{\mu}{T}\geq0\;$
. 

The input map are defined by the matrices $J_{\textrm{port}}=\left(\begin{array}{cc}
0 & g\\
-g^{\top} & 0
\end{array}\right)$ and $M_{\textrm{port}}=\left(\begin{array}{cc}
0 & g_{S}\\
-g_{S}^{\top} & 0
\end{array}\right)$ where $g^{\top}=\left(\begin{array}{cccc}
1 & 0 & 0 & 0\\
0 & 0 & 0 & 0
\end{array}\right)\;$ and $g_{S}^{\top}=\left(\begin{array}{cccc}
0 & 0 & 0 & 0\\
0 & 0 & 0 & 1
\end{array}\right).$ The extended Hamiltonian and entropy functions are: 
\begin{eqnarray*}
H_{tot}(S,V,q,p,\xi_{1},\xi_{2}) & = & H(S,V,q,p)+\xi_{1}u_{1}+\xi_{2}u_{2}\\
S_{tot}(S,V,q,p,\xi_{1},\xi_{2}) & = & S+\xi_{1}
\end{eqnarray*}
where we fixed the value $\tau_{1}=1$ associated with the port where
the irreversible phenomenon, the heat conduction with the thermostat
at temperature $u_{1}$, takes place with the positive function $\gamma_{port}=\frac{\lambda_{e}}{Tu_{1}^{2}}$
and we fixed $\tau_{2}=0$ for the port where the external force $u_{2}$
acting on the piston in a reversible way. Observe that $\{S_{tot},H_{tot}\}_{J_{port}}=u_{1}-T$
and moreover $\nabla S_{tot}$ is a Casimir of $M_{port}$.


The energy and the entropy balance equation become $\frac{dH}{dt}=y_{1}u_{1}+y_{2}u_{2}$
and $\frac{dS}{dt}-y_{1}=\frac{(u_{1}-T)^{2}}{Tu_{1}^{2}}+\frac{1}{T}\mu v^{2}\geq0\;$.

\subsection{Discrete gradient methods}

For ODEs in skew-gradient form, i.e. $\dot{x}=J(x) \frac{\partial H}{\partial x}(x)$ with $x\in \mathbb{R}^n$ and $J(x)$ a skew-symmetric matrix, we can use discretizations of the gradient $\nabla H(x)\equiv \frac{\partial H}{\partial x}$ to define a class of numerical integrators which exactly preserve the first integral $H$.
\begin{definition}\label{def31}{\cite{GONZ}}
	Let $H:\mathbb{R}^n\longrightarrow \mathbb{R}$ be a differentiable function. Then $\bar{\nabla}H:\mathbb{R}^{2n}\longrightarrow \mathbb{R}^n$ is a discrete gradient of $H$ if it is continuous and satisfies
	\begin{subequations}
		\label{discGrad}
		\begin{align}
			\bar{\nabla}H(x,x')^T(x'-x)&=H(x')-H(x)\, , \quad \, \mbox{ for all } x,x' \in\mathbb{R}^n  \, ,\label{discGradEn} \\
			\bar{\nabla}H(x,x)&=\nabla H(x)\, , \quad \quad \quad \quad \mbox{ for all } x \in\mathbb{R}^n  \, . \label{discGradCons}
		\end{align}
	\end{subequations}
\end{definition}
Some examples of discrete gradient methods can be find in  \cite{GONZ}, \cite{ITOH} and \cite{McQuRo}. 

Once a discrete gradient $\bar{\nabla}H$ has been chosen, it is straightforward to define an energy-preserving integrator by
$
\frac{x'-x}{h}=\tilde{J}(x,x',h)\bar{\nabla}H(x,x')$
where $\tilde{J}$ is a differentiable skew-symmetric matrix approximating $J$. In particular it satisfies $\tilde{J}(x,x,0)=J(x)$. 
As in the continuous case, it is immediate to check that $H$ is exactly preserved since 
$$
H(x')-H(x)=\bar{\nabla}H(x,x')^{\top}(x'-x)=h \bar{\nabla}H(x,x')^{\top} \tilde{J}(x,x',h)\bar{\nabla}H(x,x')=0 \, . 
$$
Some well-known examples of discrete gradients are:
\begin{itemize}
	\item The mean value (or averaged) discrete gradient introduced in \cite{HLvL83} and given by
	\begin{equation*}
	\label{AVF}
	\bar{\nabla}_{1}H(x,x'):=\int_0^1 \nabla H ((1-s)x+s x')\,ds 
	\end{equation*}
 for $x'\not=x$.
	\item The midpoint (or Gonzalez) discrete gradient, introduced in \cite{GONZ} and given by
	\begin{align*}
		\tiny \bar{\nabla}_{2}H(x,x')&:=\nabla H\left( \frac{1}{2}(x'+x)\right)+\tfrac{H(x')-H(x)-\nabla H\left( \frac{1}{2}(x'+x)\right)^T(x'-x)}{|x'-x|^2}(x'-x) \, , \label{gonzalez}\\
		& \mbox{ for } x'\not=x \, . \nonumber
	\end{align*}
	\item The coordinate increment discrete gradient, introduced in \cite{ITOH}, with each component given by
	\begin{equation*}
	\label{itoAbe}
	\bar{\nabla}_{3}H(x,x')_i:=\tfrac{H(x'_1,\ldots,x'_i,x_{i+1},\ldots,x_n)-H(x'_1,\ldots,x'_{i-1},x_{i},\ldots,x_n)}{x'_i-x_i}\, , \quad 1\leq i \leq N\, ,
	\end{equation*}
	when $x_i'\not=x_i$, and $\bar{\nabla}_{3}H(x,x')_i=\frac{\partial H}{\partial x_i}(x'_1,\ldots,x'_{i-1},x'_i=x_{i},x_{i+1},\ldots,x_n)$ otherwise.
\end{itemize}


\subsection{Irreversible port-Hamiltonian systems and discrete gradient methods}

Assume for simplicity that 
$J_{i}, J_{port,j}$, $M_{\alpha}$ and $M_{port,\beta}$ are constant skew symmetric matrices.  
Choose a discrete gradient $\bar{\nabla}$. Then the discretization of the IPHS  is given by 
\begin{eqnarray*}
\frac{x_{k+1}-x_k}{h}&=&
\sum_{i=1}^{k}\hat{\gamma_{i}}({x_{k+1/2}})
\{S, H\}^d_{J_{i}} J_{i}\bar{\nabla}H(x_k,x_{k+1}) + \sum_{\alpha=1}^{l} M_{\alpha}\bar{\nabla}H(x_k,x_{k+1})\\
&&
+\sum_{j=1}^{\bar{k}}\widehat{\gamma_{Port,j}}({x_{k+1/2}}, u_k)\left[ \{S_{tot}, H_{tot}\}^d_{J_{port,g_j}}\right]
g_{j} u_k + \sum_{\beta=1}^{\bar{l}}g_{S,\beta}u_{k}\\
\\
y_k&=&\sum_{j=1}^{\bar{k}}\widehat{\gamma_{Port,j}}({x_{k+1/2}}, u_k)\left[ \{S_{tot}, H_{tot}\}^d_{J_{port,g_{j}}}
\right]g_{j}^{\top}\bar{\nabla}H(x_k,x_{k+1})\\
&&+ \sum_{\beta=1}^{\bar{l}}g_{S,\beta}^{\top}\bar{\nabla}H(x_k,x_{k+1})
\end{eqnarray*}
where we denote by $\{S, H\}^d_{J_i}(x_{k}, x_{k+1})=\bar{\nabla}S(x_k,x_{k+1})^{\top} J_i\bar{\nabla}H(x_k,x_{k+1})$ and $\{S_{tot}, H_{tot}\}^d_{J_{port,g_{j}}}(x_k, x_{k+1})=(\bar{\nabla}S(x_k,x_{k+1}), \tau_k)^{\top} 
J_{port,g_{j}}
(\bar{\nabla}H(x_k,x_{k+1}), u_k)$. We also use the notation $x_{k+1/2}=(x_k+x_{k+1})/2$.
In this case, $S$ is a linear function, as a consequence,  $\bar{\nabla}S(x_k,x_{k+1})=\nabla S(x_{k})$ is in the
kernel of $M_{\alpha}$ and $\left\{ S_{\textrm{tot}},H_{\textrm{tot}}\right\}^d_{J_{\textrm{port},g_{S\beta}}}=0$, assuming we are in the conditions of Definition \ref{def:IPHS-IrrPortMap}.
Therefore, we obtain the discrete thermodynamic balance equations
$$
H(x_{k+1})-H(x_k) - hy_k^{\top} u_k =0\; ,\qquad
S(x_{k+1})-S(x_{k})-h y_k^{\top}\tau_k
\geq 0\; .
$$

\subsection{Simulating the gas-piston system using Discrete gradient methods}
As a toy example to demonstrate the properties of the discrete gradient
based integrator, we have considered the gas-piston system with the
expression of the internal energy, for a closed system (constant
number of moles $N_{0}$) 
\[
\frac{U}{U_{\textrm{ref}}}=\exp\left(\frac{1}{c}\frac{\left(S-S_{\textrm{ref}}\right)}{N_{0}R}\right)\left(\frac{V_{\textrm{ref}}}{V}\right)^{\frac{1}{c}}
\]
 where $R$ is the universal gas constant and $c$ is a positive number
depending on the nature of the gas (for a mono-atomic perfect gas
$c=\frac{3}{2}$ for instance) \cite[page 68]{Callen85}. In this
toy example we have chosen the reference values $U_{\textrm{ref}}=1$,
$V_{\textrm{ref}}=1$ and $S_{\textrm{ref}}=0$. The heat
diffusion coefficient is taken $\lambda_{e}=1$, the gravity constant
$g=10$ and the dissipation coefficient $\mu=0.5$. The initial state
of the system is chosen $S_{0}=2,V_{0}=4,q_{0}=1,p_{0}=0$. The controls
are the thermostat temperature $u_{1}=10$ and the external mechanical
force $u_{2}=-10$.

\begin{figure}
    \centering
    \begin{subfigure}[b]{0.45\textwidth}
        \centering
        \includegraphics[scale=0.45]{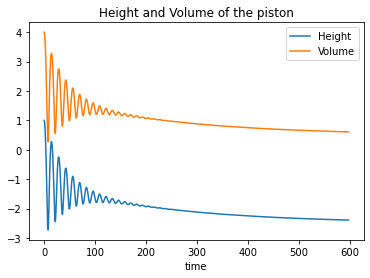}
        \caption{Height and volume vs. time}
        \label{f1}
    \end{subfigure}
    \hfill
    \begin{subfigure}[b]{0.45\textwidth}
        \centering
        \includegraphics[scale=0.45]{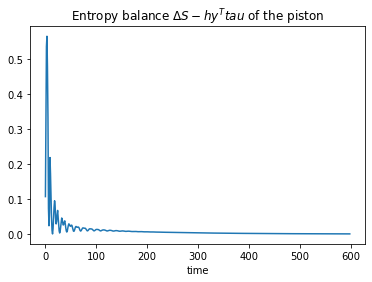}
    \caption{Positive Entropy balance in time.}
            \label{f2}
    \end{subfigure}
    \medskip
    \begin{subfigure}[b]{0.45\textwidth}
        \centering
        \includegraphics[scale=0.45]{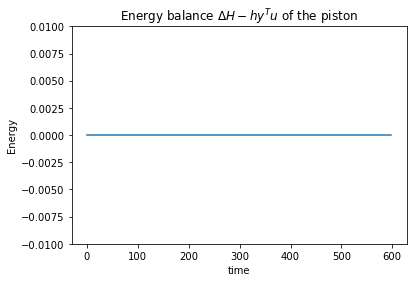}
        \caption{Energy balance vs time}
                \label{f3}
    \end{subfigure}
    \hfill
    \begin{subfigure}[b]{0.45\textwidth}
        \centering
        \includegraphics[scale=0.45]{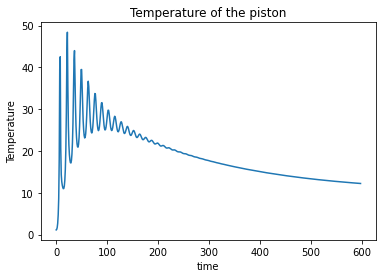}
        \caption{Piston temperature evolution}
         \label{f4}
    \end{subfigure}
    \caption{Simulation of the gas-piston system with controls.}
\end{figure}

On the figure \ref{f3}, one recovers that the numerical scheme preserves the power balance equations and on the figure \ref{f2} that the internal energy production term is indeed positive. In the figure \ref{f1} the evolution of the volume of the gas and the displacement are shown and one sees that they are proportional and can clearly distinguish two time interval, the first one, where the mechanical oscillations
is dominant, the gas acting as a nonlinear spring and the second one, corresponding to the thermal relaxation to the equilibrium. In the
figure \ref{f4}, the time evolution of the temperature is shown and the two time intervals appear clearly: the first one, corresponding
to the damped oscillatory behaviour where mechanical energy due to the external force is converted to internal energy of the gas and the second one where the thermal relaxation is dominant.


\section{Conclusions}

In this paper, we have suggested a discretization scheme adapted for irreversible port-Hamiltonian systems that allows to verify both the energy balance equation and the positivity of entropy creation terms in the entropy balance equation. In future work, we intend to extend our scheme to the conservation of other geometric properties of this class of system to more general classes of systems, e.g. systems with constraints or implicit systems and apply this scheme to more realistic physical systems.

\section*{Acknowledgements}
A. Anahory Simoes and D. Mart{\'\i}n de Diego acknowledges financial support from the Spanish Ministry of Science and Innovation, under grants PID2019-106715GB-C21 and CEX2019-000904-S.

\bibliographystyle{unsrt}
\bibliography{References}

\end{document}